\documentclass[12pt]{article}
\usepackage{t1enc}
\usepackage[latin1]{inputenc}
\usepackage[english]{babel}
\usepackage{amsmath,amsthm}
\usepackage{amsfonts}
\usepackage{latexsym}
\usepackage[dvips]{graphicx}
\usepackage{graphicx}
\usepackage[natural]{xcolor}
\usepackage{float}
\usepackage{enumerate} 
\usepackage[subnum]{cases}
\usepackage{multirow}
\usepackage{hyperref}
\usepackage[total={6in, 8.5in}]{geometry}
\newcommand{\pf}{\noindent\textbf{Proof}.\quad}
\usepackage{listings}
\makeatletter
\newcommand{\setword}[2]{%
	\phantomsection
	#1\def\@currentlabel{\unexpanded{#1}}\label{#2}%
}
\makeatother

\makeatletter
\renewcommand*\env@matrix[1][*\c@MaxMatrixCols c]{%
	\hskip -\arraycolsep
	\let\@ifnextchar\new@ifnextchar
	\array{#1}}
\makeatother

\usepackage{graphicx}
\usepackage[natural]{xcolor}
\usepackage{mathrsfs}

\usepackage{xcolor}
\long\def\ignore#1{}

\let\oldi\ignore

 \oldi{
\newtheorem{THM}{\textbf{Theorem}}[section]
\newtheorem{THMs}{\textbf{Theorem}}[section]
\newtheorem{DEF}[THM]{\textbf{Definition}}[section]
\newtheorem{LEM}[THM]{\textbf{Lemma}}
\newtheorem{CON}[THM]{\textbf{Conjecture}}
\newtheorem{PROP}[THM]{\textbf{Proposition}}
\newtheorem{COR}[THM]{\textbf{Corollary}}
\newtheorem{CORs}{\textbf{Corollary}}[section]
\newtheorem{PRO}[THM]{\textbf{Problem}}
\newcommand{\pf}{\textbf{Proof}.\quad}

\newtheorem{FAC}{\textbf{Fact}}
\newtheorem{REM}{\textbf{Remark}}
\newtheorem{OPR}{\textbf{Operation}}
\newtheorem{CLA}{\textbf{Claim}}[section]
\setcounter{CLA}{1}
 }

\newtheorem{THM}{Theorem}[section]

\newtheorem{LEM}[THM]{Lemma}

\newtheorem{CON}[THM]{Conjecture}

\newtheorem{COR}[THM]{Corollary}

\newtheorem{CLA}{Claim}[section]
\begin{document}
\title{\bf{Antimagic Orientation of Forests}}

\author{%
 Songling Shan\\
 Department of Mathematics, \\
 Illinois State Univeristy, Normal, IL 61790, USA \\
 \texttt{sshan12@ilstu.edu}
 \and
 Xiaowei Yu\thanks{Partially supported by the NSFC grants 11901252, 12031018, 11971205.}\\
  School of Mathematics and Statistics,\\
  Jiangsu Normal University,  Xuzhou,  221116, P. R. China\\
  \texttt{xwyu@jsnu.edu.cn}%
}

\date{}
\maketitle

 \begin{abstract}
An antimagic labeling of a digraph $D$ with $n$ vertices and $m$ arcs is a bijection from the set of arcs of $D$ to $\{1,2,\cdots,m\}$ such that all $n$ oriented vertex-sums are pairwise distinct, where the oriented vertex-sum of a vertex is the sum of labels of all arcs entering that vertex minus the sum of labels of all arcs leaving it.
A graph $G$ admits an antimagic orientation if $G$ has an orientation $D$ such that $D$ has an antimagic labeling.
Hefetz, M{\"{u}}tze and Schwartz conjectured every connected graph admits an antimagic orientation. In this paper, we support this conjecture by proving that any forest obtained from a given forest with at most one isolated vertex by subdividing each edge at least once admits an antimagic orientation.\\
 \smallskip
\noindent {\textbf{Key words}: Labeling; Antimagic labeling; Antimagic orientation; Forest.}
\smallskip

\noindent {Mathematics Subject Classifications: 05C78}

 \end{abstract}


\section{Introduction}

Graphs considered in this paper are simple, and digraphs in consideration are the orientations of simple graphs. Let $m$ be a positive integer and  $G$  be a graph with $m$ edges.  We use $V(G), E(G)$ to denote the vertex set and edge set of $G$, respectively.
A \emph{labeling} of a graph $G$ with $m$ edges is a bijection from $E(G)$ to the set  $\{1,2,\ldots, m\}$. 
A labeling of $G$ is \emph{antimagic} if for every two distinct vertices $u,v\in V(G)$, the sum of labels on the edges incident with $u$ differs from that of $v$.
A graph is \emph{antimagic} if it admits an antimagic labeling.

Hartsfield and Ringel in 1990 introduced antimagic labeling and they conjectured that every connected graph other than $K_{2}$ is antimagic \cite{MR1069559}, which is called the Antimagic Labeling Conjecture. They also conjectured in the same paper that every tree with at least three vertices is antimagic. The two conjectures have received intensive attention, but   remain open in general.
Some partial results on the Antimagic Labeling Conjecture can be found in~\cite{MR2096791, MR3515572, MR2478227, MR3372337,  MR3527991,  MR2174213, MR2682515,MR2510327,  MR3225299, MR4203695, MR3021347}.

Hefetz, M{\"{u}}tze and Schwartz~\cite{MR2674494} introduced a
variation of antimagic labelings, i.e., antimagic labelings on
directed graphs. An {\it antimagic labeling} of a directed graph with $m$ arcs is
a bijection from the set of arcs to the set $\{1,2,\ldots,m\}$ such that
any two oriented vertex-sums are distinct, where an {\it oriented
vertex-sum\/} of a vertex is the sum of labels of all arcs entering that vertex minus
the sum of labels of all arcs leaving it.
Given a graph $G$, we say $G$ admits an \emph{antimagic orientation} if $G$
has an orientation $D$ such that $D$ is antimagic.
Regarding antimagic orientation of graphs,
Hefetz, M{\"{u}}tze and Schwartz \cite{MR2674494} proposed the following conjecture.
\begin{CON}[\cite{MR2674494}]\label{antimagic-orientation}
	Every connected graph admits an antimagic orientation.
\end{CON}

For Conjecture~\ref{antimagic-orientation}, Hefetz, M{\"{u}}tze and Schwartz \cite{MR2674494} showed that
every orientation of a graph with order $n$ and minimum degree at least $c  \log n$ is antimagic, where $c$ is an absolute constant.
 Particularly,
they showed that every orientation of stars\,(other than $K_{1,2}$), wheels,
and complete graphs\,(other than $K_3$)
is antimagic. The conjecture is also verified for regular graphs~\cite{MR2674494,MR3879961, MR3996400},
biregular bipartite graphs~\cite{MR3723570}, Halin graphs~\cite{MR3996752}, graphs with large maximum degree~\cite{MR4141477}, graphs with a large independence set~\cite{MR4192045},  lobsters~\cite{MR4136178},
and subdivided caterpillar~\cite{ferraro2021antimagic}.  Recently, the first author~\cite{https://doi.org/10.1002/jgt.22721}
showed that the conjecture is true for all bipartite graphs with no vertex of degree 2 or 0, and for
all graphs with minimum degree at least 33.

Every bipartite antimagic  graph $G[X,Y]$ admits an antimagic orientation,
as we can direct all edges from $X$ to $Y$ and apply any of the antimagic labelings of $G$.
Thus by the result of Liang, Wong and Zhu~\cite{MR3225299} on atimagic labelings of trees, we know that every tree with at most one vertex of degree two admits an antimagic
orientation, and any tree obtained from a tree with no vertex of degree 2 by subdividing every edge  exactly once admits an antimagic
orientation.  These two results,  together with the result of the first author~\cite{https://doi.org/10.1002/jgt.22721} that every bipartite graphs with no vertex of degree 2 or 0 admits an antimagic
orientation, suggest that it is hard to find an antimagic orientation if a graph has many vertices of degree 2.
In this paper, we overcome this issue for forests
 and obtain the results below.

\begin{THM}\label{mainthm}
Let $F=(V,E)$ be a forest with at most one isolated vertex. If the set $\{v \in V\mid d_F(v) \ne 2\}$ is independent, then $F$ admits an antimagic orientation.
\end{THM}

The following result is a consequence of Theorem~\ref{mainthm}.
\begin{COR}
Let $F$ be obtained from any forest with at most one isolated vertex by subdividing each edge at least once. Then $F$ admits an antimagic orientation.
\end{COR}
\section{Notation and Preliminary Lemmas}

Let $G$ be a  graph.
 For $v\in V(G)$, $N_G(v)$ is the set of neighbors of $v$
in $G$, and
$d_G(v)=|N_G(v)|$ is the degree of $v$ in $G$.
For notational simplicity, we write $G-x$ for $G-\{x\}$.
If $F\subseteq E(G)$, then $G-F$ is obtained from $G$ by deleting all
the edges of $F$. For $F\subseteq E(\overline{G})$, $G+F$
is obtained from $G$ by adding all the edges from $F$ to $G$.
We also use $G-F$ and $G+F$ if $G$ is a digraph.  For two disjoint subsets $X,Y\subseteq V(G)$, we denote by $E_G(X,Y)$ the set of edges in $G$ with one endvertex in $X$ and the other one in $Y$ and let $e_G(X,Y)=|E_G(X,Y)|$.  A \emph{matching} $M$ in  $G$ is a set of independent edges, and we use $V(M)$ to denote the set of vertices saturated by $M$.  If $G$ is bipartite with two partite sets $X$ and $Y$, we denote $G$ by $G[X,Y]$ to emphasis the bipartitions.
For any two integers $a$ and $b$, let $[a,b]=\{ i\in \mathbb{Z} \mid   a\le i \le b\}$. If $a=1$ and $b\ge 1$, we write $[1,b]$ as $[b]$ for simplicity.


Let $G$ be a graph and $D$ be an orientation of $G$.
We denote by $A(D)$ the set of arcs of $D$. For a labeling $\tau$ on $A(D)$ and a vertex $v\in V(D)$, we use $s_{(D,\tau)}(v)$ to denote the oriented vertex-sum at $v$ in $D$ with respect to $\tau$, which is the sum of labels on all arcs entering $v$ minus the sum of labels on all arcs leaving $v$ in $D$. For simplicity, we write $s_{(D,\tau)}(v)$ as $s(v)$ if  $D$
and $\tau$ are understood.


%

We will use the lemma below to partition an integer set such that all the sum of the elements from each subset is congruent 0 modular an integer.
\begin{LEM}\label{l1}
	Let $a$ and $t$ be integers with $a\ge 0$ and $t\ge 1$, and let $r_1+r_2+\ldots +r_t$ be a partition of a positive integer $k$, where $r_i\ge 2$ for each $i\in[t]$.  Define $A=\left[1, \left\lfloor\frac{k}{2}\right\rfloor\right]\cup \left[\left\lfloor\frac{k}{2}\right\rfloor+a+1, k+a \right]$. Then $A$ can be partitioned into subsets $A_1,A_2,\ldots,A_t$  such that for every $i\in[t]$,  	$|A_i|=r_i$ and
		\begin{numcases}{}
	  \sum_{x\in A_i}x \equiv 0\pmod{k+a+1}, & if $k$ is even;  \nonumber\\
	 \sum_{x\in A_i}x \equiv 0\pmod{k+a}, & if $k$ is odd.  \nonumber
	\end{numcases}
\end{LEM}
\pf
The case for even $k$ is  Corollary 2.2 (i) from~\cite{MR3996752}. Thus we assume $k$ is odd.
Since $k=r_1+r_2+\ldots +r_t$ and $r_i\ge 2$ for each $i\in [t]$, there exists an odd $r_{i_0}\ge 3$ for some $i_0\in [t]$. Let $r'_{i_0}=r_{i_0}-1$ and $r'_i=r_i$ for $i\in [t]\backslash\{i_0\}$. Since $k-1$ is even and $k-1=r'_1+r'_2+\ldots +r'_t$,  applying the case when $k$ is even,  the set $A\backslash\{k+a\}$ can be partitioned into pairwise disjoint subsets $B_1,B_2,\ldots, B_t$ such that for every $i\in [t]$,
\[
|B_i|=r'_i\quad\quad\mbox{and}\quad\quad  \sum\limits_{x\in B_i}x \equiv 0 \pmod{k+a}.
\]
 Define
\[
A_i=\left\{
\begin{array}{ll}
	B_i, & \hbox{if $i\in [t]\backslash\{i_0\}$;} \\
	B_i\cup\{k+a\}, & \hbox{if $i=i_0$.}
\end{array}
\right.
\]
It is clear that  $|A_i|=r_i$ and $\sum\limits_{x\in A_i}x \equiv 0 \pmod{k+a}$.
\qed

This lemma below tell us how to label paths in a forest, which will be used in the proof of Theorem \ref{mainthm}.

\begin{LEM}\label{paths1}
Let $n\ge 1$ be an integer and $\overrightarrow{P_1},\overrightarrow{P_2},\ldots,\overrightarrow{P_n}$ be $n$ internally disjoint directed paths.  For each $i\in [n]$, let
$\overrightarrow{P_i}= v_{i0}v_{i1}\ldots v_{i\ell_i}$,  where $\ell_i\ge 1$ is the length of $\overrightarrow{P_i}$. Suppose  $\sum_{1\le i\le n}\ell_i=\ell$,
and the set of first edge  $v_{i0}v_{i1}$ of each $\overrightarrow{P_i}$
is labeled by a mapping $\sigma$  such that $\sigma(v_{i0}v_{i1})=i$ for every $i\in [n]$.
Then there exists 
a bijection $\tau:  \bigcup_{ 1\le i \le n} A(\overrightarrow{P_i})\rightarrow [\ell]$ satisfying the following properties:
\begin{enumerate}[(i)]
\item $\tau(v_{i0}v_{i1})=\sigma(v_{i0}v_{i1})$ for each $i\in [n]$; and
  \item $s(u)\ne s(v)$ for any two  distinct vertices  $u$ and $v$ such that $u$ is internal  of some $\overrightarrow{P_i}$
  and $v$ is internal of some $\overrightarrow{P_j}$
  with $i,j\in [n]$; and

  \item  for any  vertex $v$ that is internal of some $\overrightarrow{P_i}$ with
  $i\in [n]$, we have $| s(v)|\in [\ell-1]$ and $ s(v)\le \ell-n-1$ if  $s(v)>0$.
\end{enumerate}

\end{LEM}

\pf
Let $i\in[n]$,  $j\in [\ell_i]$, and $e_{ij}=v_{i(j-1)}v_{ij}$.  For  two distinct edges $e_{ij}$ and $e_{st}$ with $j$ and $t$ having the same parity, we write $e_{ij}\prec e_{st}$ if either $j<t$ or $j=t$ and $i<s$. Note that if $e_{ij}\prec e_{st}$, then $e_{i(j+1)}\prec e_{s(t+1)}$.  By this definition,
 all arcs from $\bigcup_{ 1\le i \le n} A(\overrightarrow{P_i})$
are ordered into two disjoint linear orderings:
\begin{eqnarray}
 L_1: e_{11} \prec e_{21} \prec  \ldots \prec{e_{n1}} \prec \ldots,   \nonumber \\
 L_2: e_{i_12} \prec e_{i_22} \prec \ldots \prec{e_{i_k2}} \prec\ldots,  \nonumber
\end{eqnarray}
where  $i_1, i_2, \ldots i_k$  with  $i_1< i_2< \ldots <i_k$ are indices such that $\overrightarrow{P_{i_j}}$ has length
at least 2 for each $j\in [k]$.
For easier assigning labels, we assume that the linear ordering $L_1$
contains $\ell_o$ edges and the linear ordering $L_2$
contains $\ell_e$ edges for some $\ell_o, \ell_e\in [\ell]$, and  we denote the $i$-th edge of $L_1$ by $f_i$,
and the $j$-th edge of $L_2$ by $g_j$ for each $i\in [\ell_o]$ and $j\in [\ell_e]$:
\begin{eqnarray}
	L_1: f_1 \prec f_2 \prec  \ldots  \prec f_{\ell_o},   \nonumber \\
	L_2: g_1 \prec g_2 \prec \ldots \prec{g_{\ell_e}}.  \nonumber
\end{eqnarray}
See Figure~1 for an illustration of such orderings.
\begin{figure}[htbp]\label{fig4}
	\begin{center}
		\scalebox{0.4}[0.4]{\includegraphics {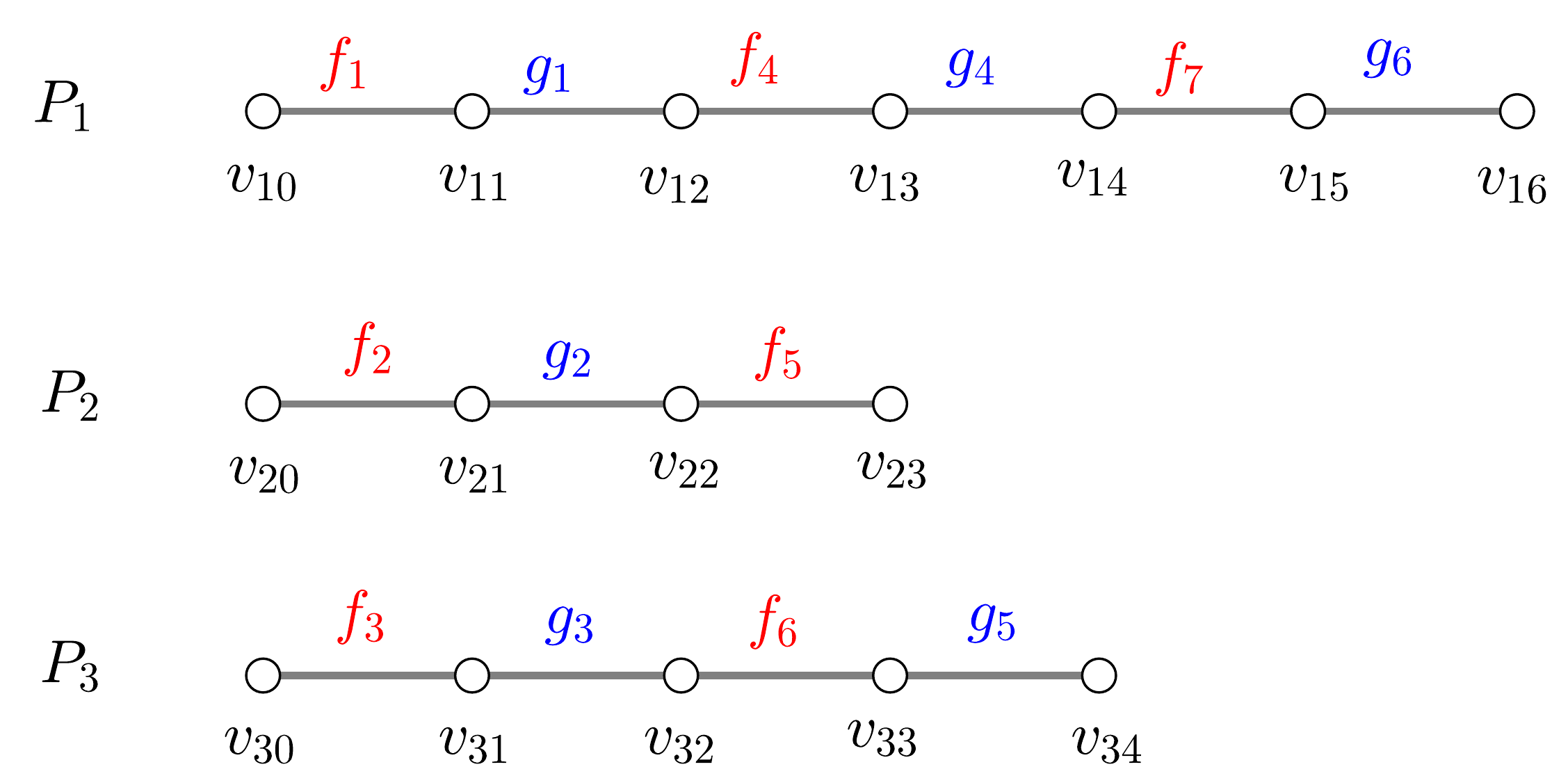}}\\
		{ Figure 1: Illustration of the orderings of edges from three directed paths}
	\end{center}
\end{figure}

%
%
%
%
%
Define $\tau:  \bigcup_{ 1\le i \le n} A(\overrightarrow{P_i})\rightarrow [\ell]$ such that
\begin{numcases}{}
\tau(f_i)=i \quad\quad \text{for $i\in[\ell_o]$};\label{a13}\\
\tau(g_i)=\ell-i+1\quad\quad \text{for $i\in[\ell_e]$}.\label{b13}
\end{numcases}
The mapping $\tau$ is a bijection. We show  below that $\tau$ satisfies properties (i) to (iii).
 By the definition of $\tau$ in~\eqref{a13}, we
have $\tau(v_{i0}v_{i1})=i=\sigma(v_{i0}v_{i1})$. Thus Lemma~\ref{paths1} (i) is true.

Since $\overrightarrow{P_{i}}$ is a directed path from $v_{i0}$ to $v_{i\ell_i}$ for any $i\in [n]$, for each $j\in [\ell_i-1]$, we have  $s(v_{ij})=\tau(e_{ij})-\tau(e_{i(j+1)})$. Therefore, by the definition of $\tau$,  $s(v_{ij})>0$ if $j$ is even and $s(v_{ij})<0$ if $j$ is odd.

We  now show Lemma~\ref{paths1} (ii).
Let $u$ and $v$ be two distinct vertices  such that $u$ is internal  of some $\overrightarrow{P_i}$
and $v$ is internal of some $\overrightarrow{P_j}$
with $i,j\in [n]$.   Assume $u=v_{is}$ and $v=v_{jt}$
for some $s\in [1,\ell_i-1]$ and $t\in [1, \ell_j-1]$.
If $s$ and $t$ have different parities, then one of $s(u)$ and $s(v)$ is positive and another is negative, and so
$s(u) \ne s(v)$. Thus we assume $s$ and $t$ have the same parity and
assume without loss of generality that $e_{is} \prec e_{jt}$.
The assumption $e_{is} \prec e_{jt}$ further implies that $e_{i(s+1)} \prec e_{j(t+1)}$.
%
Then as $\tau(e_{is})>\tau(e_{jt})$ and $\tau(e_{i(s+1)})<\tau(e_{j(t+1)})$ when $s$ is even, and $\tau(e_{is})<\tau(e_{jt})$ and $\tau(e_{i(s+1)})>\tau(e_{j(t+1)})$ when $s$ is odd, we have
\begin{eqnarray}\label{eqn2}
|s(u)|=|s(v_{is})|=|\tau(e_{is})-\tau(e_{i(s+1)})| > |\tau(e_{jt})-\tau(e_{j(t+1)})| = |s(v_{jt})|=|s(v)|.  \nonumber
\end{eqnarray}
This proves (ii).

We  lastly show Lemma~\ref{paths1} (iii). Let $v$ be an internal vertex of some $\overrightarrow{P_i}$ with
$i\in [n]$.
Since the oriented vertex-sum at $v$ is the difference of the two labels on the edges incident to $v$,
and the two  labels are distinct numbers from
 $[\ell]$,  we have $| s(v)|\in [\ell-1]$.
 Recall that $s(v)>0$  if $v=v_{ij}$ for some even $j\in [2,\ell_i-1]$.
 Since  all the labels from $[n]$ are assigned to the set of the first edges
 of the $n$ paths and $j+1 \ge 3$ is odd, we have $\tau(e_{i(j+1)}) \ge n+1$.
 Thus $s(v)=\tau(e_{ij})-\tau(e_{i(j+1)}) \le \ell -(n+1)=\ell-n-1$,
 proving the second part of (iii).
%
\qed

\section{Proof of Theorem \ref{mainthm}}
Let $m=e(F)$.  We may assume $m>0$, as otherwise,
the statement is vacuously true.
If $F$ has an isolated vertex, then its oriented vertex-sum is zero.  We will define next an orientation $D$ of $F$ and an antimagic labeling $\tau$ of $D$ such that each oriented vertex-sum at a vertex of degree at least one is non-zero. Thus, in the following, we focus only on the nontrivial components of $F$.
Let
\[
X=\{v\in V\mid d_F(v)=2 \}\quad\quad\mbox{and}\quad\quad Y=V\backslash X.
\]
Since $F$ has leaves,  we have $Y\ne \emptyset$.
Furthermore, $Y$ is independent by the condition on $F$.
 Since $Y\ne \emptyset$ and $m>0$,  we have $X\ne \emptyset$.  Because $F[X]$ is a subgraph of $F$ and all vertices in $X$ have degree 2 in $F$,  it follows that $F[X]$ is a forest with each component being a path.
 Let $P_1,P_2,\ldots,P_s$ be all the paths of $F[X]$ for some integer $s\ge 1$, and for every $i\in[s]$, let
  \[
P_i=v_{i0}v_{i1}v_{i2}\ldots v_{i\ell_i}\quad\quad \mbox{for some integer $\ell_i\ge 0$}.
\]
As $d_{P_i}(v_{i0})=d_{P_i}(v_{i\ell_i})\le  1$ and $d_{F}(v_{i0})=d_{F}(v_{i\ell_i})=2$,  it follows that each of $v_{i0}$ and  $v_{i\ell_i}$ is adjacent to a vertex  from $Y$. Let
$v_{i0}'$ and $v_{i\ell_i}'$ be the vertices from $Y$ such that
$$
v_{i0}v'_{i0}, v_{i\ell_i}v'_{i\ell_i} \in E_F(X,Y).
$$
It is possible to have  $v_{i0}=v_{i\ell_i}$, but $v'_{i0}$ and  $v'_{i\ell_i}$ are two distinct vertices as $F$ is simple and contains no cycle.

\begin{CLA}\label{claim1}
There exists a matching $M\subseteq E_F(X,Y)$ saturating   exactly one of the endvertices of $P_i$ for each $i\in [s]$ and  saturating all vertices of degree at least 3 from $Y$.
\end{CLA}
\pf We construct a new bipartite graph $H$ based on $F$. For each $P_i$ we contract it into a vertex $w_i$.
Let $V(H)=\{w_1, \ldots, w_s\} \cup Y$ and $E(H)=E_F(X,Y)$.
As the neighbor of $v_{i0}$  from $Y$ is distinct with that of $v_{i\ell_i}$, $H$
is simple such that  $d_H(w_i)=2$ for each $i\in[s]$. It suffices to show that $H$ has a matching saturating $\{w_1,\ldots, w_s\}$ and all vertices of $Y$
that have degree at least 3 in $H$. We claim first that $H$ has a matching $M_1$ saturating $\{w_1,\ldots, w_s\}$ and a matching $M_2$
saturating all vertices of $Y$
that have degree at least 3 in $H$.
Suppose $M_1$ does not exist. Then by Hall's Theorem, there is $A\subseteq \{w_1,\ldots, w_s\}$ such that
$|N_H(A)|<|A|$. This in turn implies that $e_H(A, N_H(A))=2|A| >|A|+|N_H(A)|$.
Thus $H[A,N_H(A)]$ contains a cycle, which also corresponds to a cycle of $F$,  a contradiction.
A similar argument shows that $M_2$ exists.

Let $X_1=\{w_1,\ldots, w_s\}$ and $Y_1=V(M_2)\cap Y$.
 We  claim next that we can construct a  matching $M$
that saturates both $X_1$ and $Y_1$. Let $D$
be a component of the graph formed by the union of edges from $M_1$ and $M_2$. Then $D$ can be an even cycle, a path of length even, or a path of length odd.  If $D$ is an even cycle or a path of odd length,   then $D$ has a matching $M_D$ saturating $V(D)$.     If $D$ is a path of even length and start and so end at  two vertices from $X_1$,  we let  $M_D=E(D)\cap M_1$ be a matching of $D$ (this case actually does not exist as $M_1$ saturates all vertices from $X_1$).   If $D$ is a path of even length and start and so end at  two vertices from $Y_1$,  we  let  $M_D=E(D)\cap M_2$ be a matching of $D$.   In all the three  cases, by the construction,  $M_D$ saturates all vertices in $(X_1\cap V(D)) \cup (Y_1\cap V(D))$.
Let $M$ be the union of
matchings $M_D$ for all the components $D$. By the construction, $M$ is a matching with the desired property.
\qed

By Claim~\ref{claim1}, we let $M  \subseteq E_F(X,Y)$ be a  matching saturating exactly one of the endvertices of $P_i$ for each $i\in [s]$ and all vertices of degree at least 3 from $Y$.
By the choice of $M$,  we have  $|M|=s$.
Assume
$$
Y=\{y_1, y_2, \ldots, y_{n_2}\},
$$
for some integer  $n_2 \ge 2$ ($|Y| \ge 2$ as $F$ has at least two leaves). Note that $n_2\ge |M|$ as $M\subseteq E_F(X,Y)$.
For each $y_i \not\in V(M)$, we take an  arbitrary edge $e_i \in E_F(X,Y)\setminus M$
such that $y_i$ is incident to $e_i$. Let $M^*=M\cup \{e_i\mid y_i \not\in V(M)\,\, \text{for}\,\, i\in [n_2]\}$,
 $H=F[X,Y]- M^*$, and $n_1$ be the number of vertices from $Y$ that have degree at least three in $F$.
As $F$ contains no vertex of degree 2, $Y$ has $n_2-n_1$ vertices of degree 1 in $F$.
%
%
By renaming  vertices of $Y$, assume 
\[
\left\{
  \begin{array}{ll}
    d_F(y_i)\ge 3  & \hbox{for $i\in[n_1]$;} \\
    d_F(y_i)=1  & \hbox{for $i\in[n_1+1,n_2]$.}
  \end{array}
\right.
\]
Because $d_H(y_i)=d_F(y_i)-1\ge 2$ for $i\in [n_1]$, it holds that either $e(H)=0$ or $e(H)\ge 2$.
As $e_F(X,Y)=2s$ and $|M|=s$, we have $$h:=e(H)=e_F(X,Y)-|M^*|=2s-n_2. $$

In the remainder,  we  find an orientation $D$ of $F$ and an atimagic labeling $\tau$ of $D$ in
four steps.

\begin{enumerate}[Step 1]
	\item  Orient and label $H$: direct each edge from $Y$ to $X$.
	For each $i\in [n_1]$, let $A_i$ be the set of all edges incident to $y_i$ in $H$.
	Clearly,  $|A_1|+|A_2|+\ldots+|A_{n_1}|=h$.
	Since each $y_i$ with $i\in[n_1]$ has degree at least 2 in $H$,  we have $|A_i|\ge 2$.
	By applying Lemma~\ref{l1}  with
	$t=n_1$, $a=m-s-h$ and $r_i=|A_i|$ for each $i\in[n_1]$, the set
	$$A=\left[1, \left\lfloor\frac{h}{2}\right\rfloor\right]\bigcup \left[\left\lfloor\frac{h}{2}\right\rfloor+a+1, h+a \right]=\left[1, \left\lfloor\frac{h}{2}\right\rfloor\right]\bigcup \left[m-s- \left\lceil \frac{h}{2} \right\rceil +1, m-s\right]$$
	can be partitioned into $R_1, R_2, \ldots, R_{n_1}$ such that for each $i\in [n_1]$,
	$|R_i|=|A_i|$ and $\sum_{r\in R_i} r \equiv 0 \pmod{m-s+1}$ if $h$ is even, and
	$\sum_{r\in R_i} r \equiv 0 \pmod{m-s}$ if $h$ is odd.  Label edges in $A_i$
	by integers from  $R_i$  arbitrarily such that distinct edges receive distinct colors.
	
	Let $D_1$ be the orientation of defined $H$ above and $\sigma_1$ be the labeling of $H$ defined in Step 1.
	Then for every $i\in [n_1]$, we have
	\begin{align}\label{sum1}
		s_{(D_1, \sigma_1)}(y_i)=\left\{
		\begin{array}{ll}
			-a_i\ (m-s+1), & \hbox{if $h$ is even;} \\
			-a_i\ (m-s), & \hbox{if $h$ is odd,}
		\end{array}
		\right.
	\end{align}
	for some positive integer $a_i$.

	\item Orient and label edges in $M^*\setminus M$: direct each edge from $Y$ to $X$.
	Note that $|M^*\setminus M|=n_2-s$. We assign arbitrarily the labels in $[m-\left\lceil h/2\right\rceil- n_2+1, m-\left\lceil h/2\right\rceil-s]$
	to edges in $M^*\setminus M$ such that distinct edges receive distinct labels.
	
	We let $D_2$ be the orientation of $H+ M^*\setminus M$ given through Steps 1 and 2 and $\sigma_2$ be the labeling of $D_2$ obtained through Steps 1 and 2.
	The set of labels used so far on $D_2$ is
	$$
	\left[1, \left\lfloor\frac{h}{2}\right\rfloor\right]\bigcup \left[m-n_2- \left\lceil \frac{h}{2} \right\rceil +1, m-s\right].
	$$
	By renaming the vertices of each path $P_i$, we may assume that the endvertex $v_{i\ell_i}$ of $P_i$ is saturated by $M$.
	Thus the edge $v_{i0}'v_{i0}\in E_F(X,Y)\setminus M$.
	By permuting the $s$ paths $P_1, \ldots, P_s$, we further assume
	that
	\begin{numcases}{\sigma_2(v'_{i0}v_{i0})=}
		i, & \text{ $i\in\left[1, \left\lfloor\frac{h}{2}\right\rfloor\right]$;}\label{a2} \\
		m-s+\left\lfloor h/2\right\rfloor-i+1, & \text{ $i\in\left[\left\lfloor\frac{h}{2}\right\rfloor+1,s\right]$}.\label{b2}
	\end{numcases}
	Notice that $2s= e_F(X,Y)=e(H)+|M^*|=h+n_2$.
	
	\item Orient and label edges of $F[X]$: direct each $P_i$ from $v_{i0}$ to $v_{i\ell_i}$ for $i\in [s]$, and denote the orientation by $\overrightarrow{P_i}$.
	Let $g$ be the number of paths $P_i$ with length at least one for all $i\in\left[\left\lfloor\frac{h}{2}\right\rfloor+1,s\right]$.
	Denote these $g$ paths by
	$P_{r_1}, P_{r_2}, \ldots, P_{r_g}$, where  $\left\lfloor\frac{h}{2}\right\rfloor+1 \le r_1<r_2<\ldots<r_g \le s$.
	Define a bijection  $\sigma_3^*$ from $\{v_{i_j0}v_{i_j1} \mid j\in[g]\}$ to $[\lfloor\frac{h}{2}\rfloor+1, \ \lfloor\frac{h}{2}\rfloor+g]$ such that
	\begin{align}
		\sigma^*_3(v_{r_j0}v_{r_j1})=\left\lfloor{h}/{2}\right\rfloor+j \quad\quad \text{for $j\in [g]$}. \label{d10}
	\end{align}
Denote by $\sigma_3'$  the combination of the labeling $\sigma_2$ and $\sigma_3^*$ on $D_3':=D_2+\{v_{r_j0}v_{r_j1} \mid j\in[g]\}$. Since $s_{(D'_3,\sigma'_3)}(v_{r_j0})=\sigma_2(v'_{r_j0}v_{r_j0})-\sigma^*_3(v_{r_j0}v_{r_j1})$ for $j\in[g]$,
by~\eqref{b2} and~\eqref{d10}, we have
\begin{eqnarray}\label{X2}
 m-s-(\lfloor h/2\rfloor+1) &\ge&
  s_{(D_3',\sigma_3')}(v_{r_10})>s_{(D_3',\sigma_3')}(v_{r_20})>\ldots >s_{(D_3',\sigma_3')}(v_{r_g0})   \nonumber \\
   &\ge& m-2s+\lfloor h/2\rfloor+1-(\lfloor h/2\rfloor+g) \nonumber \\
  &=& m-2s-g+1>0,
\end{eqnarray}
as  $m\ge  e_F(X,Y)+e(F[X]) \ge 2s+g$.

For $i\in \left[\left\lfloor\frac{h}{2}\right\rfloor\right]$, let ${R_i}=v'_{i0}v_{i0} {P_i}$ and $\overrightarrow{R_i}$ be the directed path
from $v'_{i0}$ to $v_{i\ell_i}$. For $i\in \left[\left\lfloor\frac{h}{2}\right\rfloor+1, \left\lfloor\frac{h}{2}\right\rfloor+g\right]$,
let $\overrightarrow{R_i}=\overrightarrow{P}_{r_{j}}$, where $j=i-\left\lfloor\frac{h}{2}\right\rfloor \in [g]$.
For every $i\in \left[\left\lfloor{h}/{2}\right\rfloor+g\right]$, we have
\[
e(R_i)\ge 1 \quad \quad\mbox{and}\quad\quad \sum\limits_{1\le j\le \left\lfloor{h}/{2}\right\rfloor+g}e(R_j)= m-|M^*|-h+\left\lfloor{h}/{2}\right\rfloor=m-n_2-\left\lceil{h}/{2}\right\rceil.
\]

We apply Lemma~\ref{paths1} on $\overrightarrow{R_1}, \ldots, \overrightarrow{R}_{\lfloor{h}/{2}\rfloor+g}$ with
\[
n:=\left\lfloor{h}/{2}\right\rfloor+g, \quad  \quad \ell:= m-n_2-\left\lceil{h}/{2}\right\rceil,\quad \text{and} \quad \sigma:=\sigma'_3
\]
to get a labeling  for those arcs of the $n$ directed paths.
Denote by $D_3$ the orientation of $F-M$ obtained through Steps 1 to 3, and let   $\sigma_3$  be the labeling of $D_3$
obtained through the three steps. Then
for any two distinct vertices  $u, v\in X\setminus \{v_{i\ell_i}, v_{r_j0}, v_{r_j\ell_{r_j}} \mid i\in [\lfloor h/2\rfloor], j\in [\lfloor h/2\rfloor+1,s]\}$,
where the set  consists of all internal vertices from the paths $\overrightarrow{R_1}, \ldots, \overrightarrow{R_n}$, by Lemma~\ref{paths1},
 we have
\begin{enumerate}[(i)]
	\item $s_{(D_3, \sigma_3)}(u)\ne s_{(D_3, \sigma_3)}(v)$; and
	\item $|s_{(D_3, \sigma_3)}(v)|\in [1, m-n_2-\left\lceil{h}/{2}\right\rceil-1]$; and
	\item if $s_{(D_3, \sigma_3)}(v)>0$, then $s_{(D_3, \sigma_3)}(v)\le m-n_2-g-h-1$.
\end{enumerate}

\item  Orient and label edges of $M$: direct each edge in $M$ from $Y$ to $X$.  Recall that $M=\{v'_{i\ell_i}v_{i\ell_i} \mid i\in[s]\}$. Without loss of generality, assume
\[
s_{({D_3},\sigma_3)}(v_{1\ell_1})\le s_{({D_3},\sigma_3)}(v_{2\ell_2})\le \ldots\le s_{({D_3},\sigma_3)}(v_{s\ell_s}).
\]
Let $\sigma_4$  be a bijection from $M$ to $[m-s+1,m]$ such that
$$\sigma_4(v'_{i\ell_i}v_{i\ell_i})=m-s+i$$
 for every $i\in[s]$.
  Denote by $D$ the orientation of $F$ obtained through Steps 1 to 4, and let   $\tau$  be the labeling of $D$
 obtained through the four steps. Then by Step 4,   for every $i\in[s]$, we have
\begin{numcases}{}
s_{({D},\tau)}(v_{i\ell_i})=s_{({D_3},\sigma_3)}(v_{i\ell_i})+\sigma_4(v'_{i\ell_i}v_{i\ell_i})\ge m-s+2.\label{a16}\\
s_{({D},\tau)}(v_{1\ell_1})< s_{({D},\tau)}(v_{2\ell_2})< \ldots< s_{({D},\tau)}(v_{s\ell_s}).\label{b16}
\end{numcases}
\end{enumerate}
Since $M$ saturates all vertices of degree at least three from $Y$ by Claim \ref{claim1},  every vertex  from  $V(M^*\setminus M)\cap Y$ is a degree 1 vertex of $F$.
By Step 4 and Step 2,
\begin{numcases}{}
	s_{(D,\tau)}(y_i)   \le -(m-s+1),   &  $ y_i\in V(M)\cap Y$;  \label{eqn10} \\
	-s_{(D,\tau)}(y_i) \in [m-\left\lceil h/2\right\rceil- n_2+1, m-\left\lceil h/2\right\rceil-s], &   \text{ $y_i\in V(M^*\setminus M)\cap Y$}.   \label{eqn11}
\end{numcases}

Next, we show that $\tau$ is an antimagic orientation of $D$. Note first that $\tau$ is a bijection from $A(D)$ to $[m]$ as  the new labels
used in each step are all distinct  and they all together form the set $[m]$, and the labelings defined in each step are all bijections.  Thus we show that for
any two distinct vertices $u,v\in V(D)$, $s_{(D,\tau)}(u) \ne s_{(D,\tau)}(v)$.
Let
$$
X_1=\{v_{i\ell_i} \mid i\in[s]\}, \quad  X_2= \{ v_{r_i0} \mid i\in [g]\}, \quad X_3=X\setminus(X_1\cup X_2).
$$
Note that  $X_1=X\cap V(M)$, and $X_3=X\setminus \{v_{i\ell_i}, v_{r_j0}, v_{r_j\ell_{r_j}} \mid i\in [\lfloor h/2\rfloor], j\in [\lfloor h/2\rfloor+1,s]\}$
is the set of internal vertices from the paths $\overrightarrow{R_1}, \ldots, \overrightarrow{R_n}$ defined in Step 3.
We have the following cases to consider:
\begin{enumerate}
	\item  $u,v\in X_1$,  $u,v\in X_2$, $u,v\in X_3$, or $u,v\in Y$;
		
	 \item  $u\in X_1$ and $v\in V(D)\setminus X_1$;

\item $u\in Y$ and $v\in X_2\cup X_3$;
	\item $u\in X_2$ and $v\in X_3$.

\end{enumerate}

{\noindent \bf Case (1)}: If $u,v \in X_1$, we get $s_{(D,\tau)}(u) \ne s_{(D,\tau)}(v)$ by~\eqref{b16}.
If $u,v\in X_2$, we have $s_{(D,\tau)}(u) \ne s_{(D,\tau)}(v)$ by~\eqref{X2}.
If $u, v\in X_3$, we have
$s_{(D,\tau)}(u) \ne s_{(D,\tau)}(v)$ by Property (i) of $\sigma_3$ defined in Step 3.
Thus we assume $u,v\in Y$. If $u,v\in Y\setminus V(M)$, then $u$ and $v$ are leaves of $F$ by the choice of $M$
and so $s_{(D,\tau)}(u) \ne s_{(D,\tau)}(v)$. If $u\in Y\setminus V(M)$ and $v\in Y\cap V(M)$, then $s_{(D,\tau)}(u) > s_{(D,\tau)}(v)$
by~\eqref{eqn10} and~\eqref{eqn11}. Hence we assume $u,v\in Y\cap V(M)$. Then by~\eqref{sum1} and~\eqref{eqn10}, we have
	\begin{align} \nonumber
	s_{(D,\tau)}(u) - s_{(D,\tau)}(v) \equiv \left\{
	\begin{array}{ll}
		b \pmod{m-s+1}, & \hbox{if $h$ is even;} \\
	b \pmod{m-s}, & \hbox{if $h$ is odd};
	\end{array}
	\right.
\end{align}
where $b$ is an integer satisfying  $b\in [-(s-1), s-1]$ and $b\ne 0$. Since $m\ge 2s$,  we have $b<m-s$.
Thus $s_{(D,\tau)}(u) \ne s_{(D,\tau)}(v)$.

\smallskip

{\noindent \bf Case (2)}: Let $u\in X_1$ and $v\in V(D)\setminus X_1$. By~\eqref{a16}, we have $s_{(D,\tau)}(u) \ge m-s+2$.
If $v\in Y$, we have $s_{(D,\tau)}(v) <0$. Thus we assume $v\in X_2\cup X_3$.
By~\eqref{X2} and Property (ii) of $\sigma_3$ defined in Step 3,  we have
$s_{(D,\tau)}(v) \le \max\{m-s-(\lfloor h/2\rfloor+1), m-n_2-\left\lceil{h}/{2}\right\rceil-1\}<m-s $, as $n_2\ge s$.
Thus $s_{(D,\tau)}(u) \ne s_{(D,\tau)}(v)$.
	
\smallskip

{\noindent \bf Case (3)}: Let	$u\in Y$ and $v\in X_2\cup X_3$.
Note that $s_{(D,\tau)}(u) \ne s_{(D,\tau)}(v)$ if $v\in X_2$, as $s_{(D,\tau)}(u)<0$ and $s_{(D,\tau)}(v)>0$.
Thus we assume $v\in X_3$.
By~\eqref{eqn10} and~\eqref{eqn11}, we have $|s_{(D,\tau)}(u)| \ge m-n_2-\lceil{h}/{2}\rceil+1$. By Property (ii) of $\sigma_3$ defined in Step 3,  we have
$|s_{(D,\tau)}(v)| \le m-n_2-\left\lceil{h}/{2}\right\rceil-1 < |s_{(D,\tau)}(u)|$. Thus $s_{(D,\tau)}(u) \ne s_{(D,\tau)}(v)$.

\smallskip

{\noindent \bf Case (4)}: Let	$u\in X_2$ and $v\in X_3$. By~\eqref{X2}, we have $s_{(D,\tau)}(u) \ge m+1-2s-g-h>0$.
Thus we only consider  $v\in X_3$ such that $s_{(D,\tau)}(v)>0$.
By Property (iii) of $\sigma_3$ defined in Step 3, we have  $s_{(D, \tau)}(v)\le m-n_2-g-h-1=m-2s-g-1<m+1-2s-g$, where recall $2s=h+n_2$.
Thus   $s_{(D,\tau)}(u) \ne s_{(D,\tau)}(v)$.

The proof is complete. \qed

\bibliographystyle{abbrv}

\bibliography{Bibifile}

\end{document}